\documentclass[final, a4paper, 10pt]{article}

\usepackage{theorem}
\usepackage{amsfonts}
\usepackage{amsbsy}
\usepackage{amsmath}
\usepackage{latexsym}
\usepackage{verbatim}
\usepackage{graphicx}
\usepackage{pb-diagram}
\newtheorem{thm}{Theorem}[section]
\newtheorem{cor}[thm]{Corollary}
\newtheorem{lem}[thm]{Lemma}
\newtheorem{prop}[thm]{Proposition}
\theorembodyfont{\rmfamily}
\newtheorem{rem}[thm]{Remark}
\theorembodyfont{\rmfamily}
\newtheorem{exmp}[thm]{Example}
\numberwithin{equation}{section}

\newenvironment{prf}[1][Proof]{\textbf{#1.} } {$\square$ \bigskip}
\newcommand{\MTRX} {\mathbf{M}}                         %
\newcommand{\IC} {\mathbb{C}}                           %
\newcommand{\IR} {\mathbb{R}}                           %
\newcommand{\IT} {\mathbb{T}}                           %
\newcommand{\IZ} {\mathbb{Z}}                           %
\newcommand{\BDD}[1]{\mathbf B(#1)}                     %
\newcommand{\CPT}[1]{\mathbf K(#1)}                     %
\newcommand{\UTWO} {\mathbf U (2)}                      %
\newcommand{\CA} {\mathcal{A}}                          %
\newcommand{\CD} {\mathcal{D}}                          %
\newcommand{\CK} {\mathcal{K}}                          %
\newcommand{\CL} {\mathcal{L}}                          %
\newcommand{\CM} {\mathcal{M}}                          %
\newcommand{\CS} {\mathcal{S}}                          %
\newcommand{\CT} {\mathcal{T}}                          %
\newcommand{\tr} {\mathrm{tr}}                          %
\newcommand{\Tr} {\mathrm{Tr}}                          %
\newcommand{\ind} {\mathrm{index}}                      %
\newcommand{\wn} {\mathrm{wn}}                          %
\newcommand{\df} {\mathrm{d}}                           %
\newcommand{\Inv} {\mathrm{Inv}}                        %
\newcommand{\Com} {\mathrm{Com}}                        %
\newcommand{\cstr} {\mathrm{C}^*}                       %
\newcommand{\To} {\longrightarrow}                      %
\newcommand{\ot} {\otimes}                              %
\newcommand{\by}{{\times}}                              %
\newcommand{\id} {\iota}                                %
\newcommand{\la} {\lambda}                              %
\newcommand{\vfi} {\varphi}                             %
\newcommand{\om} {\omega}                               %
\newcommand{\OM} {\Omega}                               %
\newcommand{\set}[1]{\left\{#1\right\}}                 %
\begin{document}

\title{An Index Theorem for Toeplitz Operators on the Quarter-Plane}

\author{Adel B. Badi\footnote{Mathematics Department, Faculty of Science, The 7th of October University, Misurata, Libya.\newline email: adbabadi@yahoo.com
}}

%___________________________________________________________________________________________________________________________

\maketitle

%___________________________________________________________________________________________________________________________
%

\begin{abstract}
We prove an index theorem for Toeplitz operators on the quarter-plane using the index theory for generalized Toeplitz operators introduced by G. J. Murphy in \cite{MUR1}. To prove this index theorem we construct an indicial triple on the tensor product of two $\cstr$-algebras provided with indicial triples with general conditions. We show that our results can be extended to some extensions of the theory of Toeplitz operators on the quarter-plane.
\end{abstract}

%___________________________________________________________________________________________________________________________
%

\section{Introduction}
\label{sec1}
In this section we introduce some problems emerged in \cite{DOHO1} when the classical theory of Fredholm index was used to study the index theory of Toeplitz operators on the quarter-plane. We start with setting some notation that we are going to use throughout this paper. We will denote the fields of real and complex numbers by $\IR$ and $\IC$, respectively, and the set of integers will be denoted by $\IZ$. The unit circle in $\IC$ will be denoted by $\IT$ and the function $\IT \to \IC, \ \la \mapsto \la$ will always be denoted by $z$. If $K$ is a compact space, then we will denote the algebra of continuous complex-valued functions on $K$ by $C(K)$.

Let $L^2(\IT)$ be the space of square-integrable functions relative to the normalized Lebesgue measure on $\IT$. Then the set $\set{e_m \mid m \in \IZ}$ is an orthonormal basis for the Hilbert space $L^2(\IT)$, where $e_m = z^m$. The classical Hardy space $H^2(\IT)$ is defined as the closed linear subspace of $L^2(\IT)$ spanned by the set $\set{e_m \mid m \geq 0}$. The orthogonal projection of $L^2(\IT)$ onto $H^2(\IT)$ will be denoted by $P_\IT$.

Similarly, we define $L^2(\IT^2)$ as the space of square-integrable functions relative to the normalized Lebesgue measure on $\IT^2$. We identify $L^2(\IT^2)$ with $L^2(\IT) \ot L^2(\IT)$. Hence, the set $\set{e_m \ot e_n \mid m, n \in \IZ}$ is an orthonormal basis for the Hilbert space $L^2(\IT^2)$. We define the Hardy space $H^2(\IT^2)$ as the closed linear subspace of $L^2(\IT^2)$ spanned by $\set{e_m \ot e_n \mid m, n \geq 0}$. Thus, $H^2(\IT^2) = H^2(\IT) \ot H^2(\IT)$. %Note that the dual group of the compact abelian group $\IT^2$ is the group $\IZ^2$. The space $H^2(\IT^2)$ consists of the functions in $L^2(\IT^2)$ whose Fourier transform vanishes off the semigroup $\set{(m, n) \mid m, n \geq 0}$.
Let $P_{\IT^2}$ be the orthogonal projection of $L^2(\IT^2)$ onto $H^2(\IT^2)$. Then $P_{\IT^2} = P_\IT \ot P_\IT$.

Now we define Toeplitz operators on the quarter-plane. Let $\vfi \in C(\IT^2)$. Then the Toeplitz operator $W_\vfi$ on $H^2(\IT^2)$ with symbol $\vfi$ is defined by the relation $W_\vfi(f) = P_{\IT^2}(\vfi f) \ \ (f \in H^2(\IT^2))$. These operators and the $\cstr$-algebra $\CT_{\IT^2}$ generated by them were studied by R. G. Douglas and R. Howe among others. Douglas and Howe showed that there are Fredholm operators in $\CT_{\IT^2}$ of arbitrary integer index \cite[Corollary of Theorem 1]{DOHO1}. However, they found that the Toeplitz operator $W_\vfi$ is a Fredholm operator if and only if $\vfi$ is non-vanishing and homotopic (in the set of invertible functions in $C(\IT^2)$) to a constant function
\cite[Corollary of Theorem 1]{DOHO1}. Equivalently, $\vfi = e^\psi$ for some $\psi \in C(\IT^2)$. In this case the Fredholm index is zero. Thus, the classical Fredholm theory gives a trivial index for Toeplitz operators on the quarter-plane. This situation is similar to that in \cite{BECO1} where C. A. Berger and L. A. Coburn studied Toeplitz operators on the group $\UTWO$ of unitary $2 \by 2$ matrices. In \cite{MUR2}, G. J. Murphy introduced a satisfactory index theorem for Toeplitz operators on the group $\UTWO$ using a generalization of Fredholm index he developed earlier in \cite{MUR4} and later he put his results in a general framework \cite{MUR1}. Murphy's results follow from the existence of a structure he called an indicial triple. Thus, the index theorems will follow once we have constructed a suitable indicial triple. An indicial triple on a $\cstr$-algebra defines a topological index on the given $\cstr$-algebra. Moreover, it defines a generalized Fredholm index on the (generalized) Toeplitz operators defined by the given $\cstr$-algebra. The relation between both of these indices is similar to the relation between the topological and the Fredholm indices in classical Toeplitz operators (Toeplitz operators on the unit circle $\IT$).

In this paper, we use the results in \cite{MUR1} to obtain an index theorem for Toeplitz operators on the quarter-plane. To do this we will prove in Section \ref{sec2} a general theorem showing that we can construct an indicial triple on the tensor product of two $\cstr$-algebras provided with indicial triples satisfying certain conditions. In Section \ref{sec3}, we use our construction to prove an index theorem for Toeplitz operators on the quarter-plane. In our index theorem, a Toeplitz operators $W_\vfi$ is a Fredholm operator (in a generalized sense) if and only if $\vfi$ is non-vanishing on $\IT^2$. Moreover, the Fredholm index of $W_\vfi$ is non-zero unless $\vfi = e^\psi$ for some $\psi \in C(\IT^2)$.

Now we describe some extensions of the theory of Toeplitz operators on the quarter-plane. Let $L^2(\IT^N)$ be defined in the same manner as $L^2(\IT^2)$ and let $H^2(\IT^N)$ be defined analogously, where $N > 0$. Hence, $L^2(\IT^N) = \bigotimes_{k = 1}^N L^2(\IT)$ and $H^2(\IT^N) = \bigotimes_{k = 1}^N H^2(\IT)$. We define the projection $P_{\IT^N}$ as the orthogonal projection of $L^2(\IT^N)$ onto $H^2(\IT^N)$. Thus, $P_{\IT^N} = P_\IT \ot \cdots \ot P_\IT$. The Toeplitz operator $W_\vfi$ on $H^2(\IT^N)$ with symbol $\vfi \in C(\IT^N)$ is defined by the relation $W_\vfi (f) = P_{\IT^N} (\vfi f)$, for any $f \in H^2(\IT^N)$. Douglas and Howe showed that their results on the quarter-plane do not extend to Toeplitz operators on $H^2(\IT^N)$ for $N > 2$ \cite[p. 215]{DOHO1}. Our results can be extended quite easily to Toeplitz operators on $H^2(\IT^N)$ for any positive integer $N$.

Let $L_N^2(\IT^2) = \bigoplus_{k = 1}^N L^2(\IT^2)$ and let $H_N^2(\IT^2)$ be defined similarly, where $N$ is a positive integer. Then the orthogonal projection $P_N$ of $L_N^2(\IT^2)$ onto $H_N^2(\IT^2)$ is equal to $P_{\IT^2} \oplus \cdots \oplus P_{\IT^2}$. Now we define another extension of Toeplitz operators on the quarter-plane. These operators are defined by the relation $W_\Phi(f) = P_N(\Phi f), \ f \in H_N^2(\IT^2)$, where $\Phi = [\vfi_{ij}], \ \vfi_{ij} \in C(\IT^2), \ 1 \leq i, j \leq N$. The results in \cite{DOHO1} established for Toeplitz operators on the quarter-plane do not extend to these operators if $N >1$ \cite[p. 214]{DOHO1}. On the other hand, our results can be extended with a little effort to these operators using our results on Toeplitz operators on the quarter-plane and results from~\cite{MUR1}.

%___________________________________________________________________________________________________________________________
%

\section{Index Theory for Tensor Products}
\label{sec2}
Suppose that $A$ is a $\cstr$-algebra and $\tr$ is a trace on $A$. We will denote the ideal of definition of $\tr$ by $\CM_\tr$ and the closure of $\CM_\tr$ by $\CK_\tr$. If $a, b \in A$, then we say that $[a, b] = ab - ba$ is a {\em commutator} of $a, b$. The closed ideal generated by all commutators in $A$ is called the {\em closed commutator ideal} of $A$ and we will denote it by $\Com(A)$. We denote the group of invertible elements in $A$ by $\Inv(A)$. All tensor products of $\cstr$-algebras in this paper are the spatial tensor products. The symbol~$\odot$ denotes the algebraic tensor product.

To introduce the main problem in this paper, we give the definitions and the results of the index theory for generalized Toeplitz operators developed by Murphy in \cite{MUR1}. Let $C$ and $\CL$ be unital $\cstr$-algebras such that $\CL$ contains $C$ as a unital $\cstr$-subalgebra. Suppose that $F$ is a self-adjoint unitary in $\CL$ and $\tr$ is a lower semicontinuous trace on $\CL$. If the unital $*$-subalgebra $\set{\vfi \in C \mid [F, \vfi] \in \CM_\tr}$ is dense in $C$, then we call the triple $\OM = (\CL, F, \tr)$ an {\em indicial triple} for $C$~\cite[p. 266]{MUR1}. We usually denote the latter dense $*$-subalgebra by $C_\OM$, and the commutator $[F, \vfi]$ will be denoted by $\df \vfi$, where $\vfi \in C$. Sometimes we refer to $C$ by the {\em symbol algebra} of $\OM$.

A {\em topological index} on a unital $\cstr$-algebra $C$ is a
locally constant mapping $\om : \Inv(C) \to \IR$ such that $\om(\vfi \psi) = \om(\vfi) + \om(\psi)$, for all $\vfi, \psi \in C$. Let $C$ be a $\cstr$-algebra and let $\OM = (\CL, F, \tr)$ be an indicial triple for $C$. Then there exists a unique topological index $\om$ on $C$ such that $\om(\vfi) = \frac{1}{2} \tr(\vfi^{-1} \df \vfi)$, for all $\vfi \in \Inv(C_\OM)$ \cite[Theorem 2.6]{MUR1}. %This is the most important property of indicial triples.

Let $\CL_\IT$ be the unital $\cstr$-subalgebra of $\BDD{L^2(\IT)}$ generated by $C(\IT)$ and the projection $P_\IT$. If we define the trace $\tr_\IT$ as the restriction of the canonical trace on $\BDD{L^2(\IT)}$ to $\CL_\IT$ and we define $F_\IT$ as the self-adjoint unitary $2P_\IT - 1$, then, by \cite[Example 2.7]{MUR1}, $\OM_\IT = (\CL_\IT, F_\IT, \tr_\IT)$ is an indicial triple for the $\cstr$-algebra $C(\IT)$. Let $\om_\IT$ be the topological index related to the indicial triple $\OM_\IT$. If $\vfi \in \Inv(C(\IT))$, then $\om_\IT(\vfi) = \wn(\vfi)$, where $\wn(\vfi)$ is the winding number of $\vfi$.

Let $\OM = (\CL, F, \tr)$ be an indicial triple for the $\cstr$-algebra $C$ and let $P = (F + 1)/2$. If $\vfi \in C$, then we call $T_\vfi = P \vfi P$ a {\em Toeplitz element} associated to $\OM$. The $\cstr$-algebra $\CA$ generated by all these elements is called the {\em Toeplitz algebra} associated to the indicial triple $\OM$. Note that $\CA$ is a unital $\cstr$-algebra with $P$ as its unit. We define the trace $\Tr$ on $\CA$ as the restriction of $\tr$ to $\CA$. To avoid trivialities, the trace $\Tr$ must be non-finite and this is equivalent to the condition $P \notin \CM_\tr$. Henceforth, we will assume that any indicial triple in this paper satisfies the latter condition. An element $a \in \CA$ is called an {\em $\OM$-Fredholm} element if there exists an element $b \in \CA$ such that $P - ab, P - ba \in \CM_\Tr$. The element $b$ is called a {\em partial inverse} of~$a$. The {\em $\Tr$-index} of $a$, or the {\em Fredholm index} of $a$ relative to $\Tr$, is defined by the relation $\ind_\OM(a) = \Tr(ab - ba)$. Moreover,  $a$ is an $\OM$-Fredholm element if and only if there exists an element $c \in \CA$ such that $P - ac, P - ca \in \CK_\Tr$. By \cite[Theorem 3.1]{MUR1}, if $\vfi$ is an invertible element of $C$, then $T_\vfi$ is an $\OM$-Fredholm element and $\ind_\OM(T_\vfi) = - \om(\vfi)$, where $\om$ is the unique topological index related to $\OM$.

We will refer to the $\cstr$-algebra $\CT_\IT$ generated by Toeplitz operators on the unit circle by the {\em classical Toeplitz algebra}. This algebra is $*$-isomorphic to the Toeplitz algebra $\CA_\IT$ associated to the indicial triple $\OM_\IT$ defined above. Hence, we can identify these $\cstr$-algebras. If $\vfi \in C(\IT)$, then the Toeplitz element $T_\vfi$ is an $\OM_\IT$-Fredholm element if and only if it is a Fredholm operator in the classical sense. Moreover, in this case the $\Tr_\IT$-index is identical to the classical Fredholm index, where $\Tr_\IT$ is the restriction of the trace $\tr_\IT$ to $\CA_\IT$.

Similar to the case of classical Toeplitz operators and Toeplitz operators on $\UTWO$, we want to construct an indicial triple for $C(\IT^2) = C(\IT) \ot C(\IT)$ to obtain an index theorem for Toeplitz operators on the quarter-plane. We can view $C(\IT^2)$ as a unital $\cstr$-subalgebra of $\BDD{L^2(\IT^2)}$. Let $F_{\IT^2}$ be the self-adjoint unitary $2P_{\IT^2} - 1$ and let $\CL_{\IT^2} = \CL_\IT \ot \CL_\IT \subseteq \BDD{L^2(\IT^2)}$. Note that $C(\IT^2) \subseteq \CL_{\IT^2}$ and $P_{\IT^2} \in \CL_{\IT^2}$. To~construct an indicial triple for $C(\IT^2)$, we need to define a trace on $\CL_{\IT^2}$. The~first thing comes to mind is the construction $\tr_\IT \ot \tr_\IT$. Let $\tr_0$ be the restriction of the canonical trace on $\BDD{L^2(\IT^2)}$ to $\CL_{\IT^2}$. Then $\CM_{\tr_\IT} \odot \CM_{\tr_\IT} \subseteq \CM_{\tr_0}$, $\CK_{\tr_\IT} \ot \CK_{\tr_\IT} = \CK_{\tr_0}$ and
\[
\tr_0(a \ot b) = \tr_\IT(a) \tr_\IT(b) \quad (a, b \in \CM_{\tr_\IT}).
\]
Moreover, $\tr_0$ is the only lower semicontinuous trace satisfying these conditions. Hence, we can consider the trace $\tr_0$ as the tensor product of the trace $\tr_\IT$ by itself. The trace $\tr_0$ with the $\cstr$-algebra $\CL_{\IT^2}$ and the self-adjoint unitary $F_{\IT^2}$ does not define an indicial triple for $C(\IT^2)$. To see this, suppose that $(\CL_{\IT^2}, F_{\IT^2}, \tr_0)$ is an indicial triple for $C(\IT^2)$. Then $\set{\vfi \in C(\IT^2) \mid \df \vfi \in \CM_{\tr_0}}$ is dense in $C(\IT^2)$. Consider the function $z \ot 1 \in C(\IT^2)$. There is a sequence $(\vfi_m)_{m = 1}^\infty$ in $C(\IT^2)$ such that $\df \vfi_m \in \CM_{\tr_0}$ and $\vfi_m \To z \ot 1$. This implies
\[
2[P_\IT, z] \ot P_\IT = \df (z \ot 1) \in \CK_{\tr_0} \subseteq \CPT{L^2(\IT^2)},
\]
a contradiction.
%As we mentioned above, in this paper we will use the results in \cite{MUR1} to introduce an index theorem for Toeplitz operators on the quarter-plane.

It is obvious now that the main difficulty in the construction of an indicial triple for $C(\IT^2)$ is finding a suitable trace. Similarly, the construction of the trace for defining a generalized Fredholm index for Toeplitz operators on $\UTWO$ is the main result in \cite{MUR2}. If we have a character $\tau$ on $\CL_\IT$ (a non-zero multiplicative linear functional), then it is easy to show that $\tau \ot \tr_\IT$ and $\tr_\IT \ot \tau$ are traces on $\CL_\IT \ot \CL_\IT$. Moreover, if the character $\tau$ does not vanish at $P_\IT$, then, as we will see soon, the $\cstr$-algebra $\CL_{\IT^2}$ and the self-adjoint unitary $F_{\IT^2}$ provided with the trace $\tau \ot \tr_\IT + \theta \tr_\IT \ot \tau$ define an indicial triple for $C(\IT^2)$, where $\theta$ is a fixed positive number. We will prove our result in a general setting to make it useful in other possible extensions of the theory.

Let $\OM_i = (\CL_i, F_i, \tr_i)$ be an indicial triple for the $\cstr$-algebra $C_i$ and $i = 1, 2$. We denote by $C_{\OM_i}$ the dense unital $*$-subalgebra of $C_i$ consisting of the elements $\vfi \in C_i$ such that $\df \vfi = [F_i, \vfi] \in \CM_{\tr_i}$. The main target in this section is the construction of an indicial triple $\OM = (\CL, F, \tr)$ for the $\cstr$-algebra $C = C_1 \ot C_2$. Motivated by the previous discussion, we let $\CL = \CL_1 \ot \CL_2$ and we define $F \in \CL$ as the self-adjoint unitary $2 (P_1 \ot P_2) - 1$, where ${P_i = (F_i + 1)/2}$, $i = 1, 2$. Let $\tau_i$ be a character on $\CL_i$ such that $\tau_i(P_i) \neq 0$ (equivalently, $\tau_i(P_i) = 1$), for $i = 1, 2$. Then $\tr_1 \ot \tau_2$ is a lower semicontinuous trace on $\CL_1 \ot \CL_2$, where $(\tr_1 \ot \tau_2)(a) = \tr_1((\id \ot \tau_2)(a)) \ (a \in \CL^+)$. The trace $\tau_1 \ot \tr_2$ is defined similarly. Hence, the mapping
\[
\tr: \CL^+ \to [0, +\infty], \ \tr(a) = \theta_1 (\tr_1 \ot \tau_2) (a) + \theta_2 (\tau_1 \ot \tr_2) (a) \ \ (a \in \CL^+)
\]
is a lower semicontinuous trace on $\CL$, where $\theta_1, \theta_2$ are fixed positive real numbers. We need the following lemma to show that $\CL$ and $F$ provided with this trace define an indicial triple for $C$.

\begin{lem} \label{lbl4001}
Suppose that $A, B$ are $\cstr$-algebras. Let $\tr_A, \tr_B$ be lower semicontinuous traces and let $\tau_A, \tau_B$ be characters on $A, B$, respectively. If $a \in A$ and $b \in \CM_{\tr_B} \cap \ker(\tau_B)$, then ${a \ot b}$ belongs to the definition ideal of the trace $\theta_1 \tr_A \ot \tau_B + \theta_2 \tau_A \ot \tr_B$, where $\theta_1, \theta_2$ are fixed positive numbers. Moreover,
\[
(\tr_A \ot \tau_B) (a \ot b) = 0.
\]
\end{lem}

\begin{prf}
It is obvious that $|b|$ belongs to the $\cstr$-subalgebra $\ker(\tau_B)$ of~$B$. Note that $|a \ot b| = |a| \ot |b|$. Hence,
{\setlength\arraycolsep{2pt}
\begin{eqnarray*}
   (\tr_A \ot \tau_B) \big( |a \ot b| \big) &=& \tr_A \Big( (\id \ot \tau_B) \big( |a| \ot |b| \big) \Big) \\
                                            &=& \tr_A \Big( \tau_B \big( |b| \big) |a| \Big) \\
                                            &=& \tr_A \big( 0 \cdot |a| \big) = 0.
\end{eqnarray*}}%
Thus, $|a \ot b| \in \CM_{\tr_A \ot \tau_B}^+$. By \cite[Proposition A1(b)]{PHRA1}, $\big| (\tr_A \ot \tau_B) (a \ot b) \big| \leq (\tr_A \ot \tau_B) \big( |a \ot b| \big)$. Hence, $(\tr_A \ot \tau_B) (a \ot b) = 0$. By \cite[Proposition A1(a)]{PHRA1}, $|b| \in \CM_{\tr_B}$ and it follows immediately that $|a \ot b| = |a| \ot |b| \in \CM_{\tau_A \ot \tr_B}^+$. Thus, $|a \ot b|$ belongs to the definition ideal of the trace $\theta_1 \tr_A \ot \tau_B + \theta_2 \tau_A \ot \tr_B$ and again by \cite[Proposition A1(a)]{PHRA1}, $a \ot b$ belongs to this ideal too.
\end{prf}

Now, we show that $\df \vfi \in \CM_\tr$, for every $\vfi$ in the unital dense $*$-subalgebra ${C_{\OM_1} \odot C_{\OM_2}}$ of $C$. To do this, it is enough to show that $\df (\vfi_1 \ot \vfi_2) \in \CM_\tr$, for all $\vfi_1 \in C_{\OM_1}$, $\vfi_2 \in C_{\OM_2}$. It is easy to see that
\[
\df (\vfi_1 \ot \vfi_2) = 2 \big( [P_1, \vfi_1] \ot P_2 \vfi_2 + \vfi_1 P_1 \ot [P_2, \vfi_2] \big).
\]
By Lemma \ref{lbl4001}, $[P_1, \vfi_1] \ot P_2 \vfi_2$ and $\vfi_1 P_1 \ot [P_2, \vfi_2]$ belong to the definition ideal of the trace $\tr = \theta_1 \tr_1 \ot \tau_2 + \theta_2 \tau_1 \ot \tr_2$. Hence, $\df (\vfi_1 \ot \vfi_2) \in \CM_\tr$. It follows immediately now that $\OM = (\CL, F, \tr)$ is an indicial triple for $C$.

The following theorem describes the construction of the trace on the tensor products and it gives a simple formula to calculate the index for elementary tensors of invertible elements.

\begin{thm} \label{lbl4002}
Let $\OM_i = (\CL_i, F_i, \tr_i)$ be an indicial triple for the $\cstr$-algebra $C_i$ and let $P_i = (F_i + 1)/2$, where $i = 1, 2$. If there exists a character $\tau_i$ on $\CL_i$ such that ${\tau_i(P_i) \neq 0}$, then $\OM = (\CL, F, \tr)$ is an indicial triple for $C = C_1 \ot C_2$, where $\CL = \CL_1 \ot \CL_2$, $F = 2(P_1 \ot P_2) - 1$ and the trace $\tr : \CL^+ \to [0, +\infty]$ is defined as
\[
\tr(a) = \theta_1 (\tr_1 \ot \tau_2) (a) + \theta_2 (\tau_1 \ot \tr_2) (a) \ \ (a \in \CL^+),
\]
where $\theta_1, \theta_2$ are fixed positive numbers. Moreover,
\[
 \om(\vfi_1 \ot \vfi_2) = \theta_1 \om_1(\vfi_1) + \theta_2 \om_2(\vfi_2),
\]
for all $\vfi_1 \in \Inv(C_1)$, $\vfi_2 \in \Inv(C_2)$, where $\om_1, \om_2, \om$ are the topological indices related to the indicial triples $\OM_1, \OM_2, \OM$, respectively.
\end{thm}

\begin{prf}
By the discussion above, the only thing we need to prove is the last equality, and by continuity, it is sufficient to show that
\[
 \om(\vfi_1 \ot \vfi_2) = \theta_1 \om_1(\vfi_1) + \theta_2 \om_2(\vfi_2),
\]
for all $\vfi_1 \in \Inv(C_{\OM_1})$, $\vfi_2 \in \Inv(C_{\OM_2})$. Note that%
{\setlength\arraycolsep{2pt}
\begin{eqnarray*}
  \om(\vfi_1 \ot \vfi_2) &=& \frac{1}{2} \tr \big( (\vfi_1 \ot \vfi_2)^{-1} \df (\vfi_1 \ot \vfi_2) \big) \\
   &=& \tr \big( \vfi_1^{-1}[P_1, \vfi_1] \ot \vfi_2^{-1} P_2 \vfi_2 + P_1 \ot \vfi_2^{-1} [P_2, \vfi_2] \big).
\end{eqnarray*}}%
We calculate each item separately. It is obvious that $\vfi_1^{-1} [P_1, \vfi_1] \in \ker(\tau_1)$ and $\vfi_2^{-1}[P_2, \vfi_2]  \in \ker(\tau_2)$. Hence, by Lemma \ref{lbl4001}, we have
\[
(\tr_1 \ot \tau_2) \big( P_1 \ot \vfi_2^{-1} [P_2, \vfi_2] \big) = (\tau_1 \ot \tr_2) \big( \vfi_1^{-1}[P_1, \vfi_1] \ot \vfi_2^{-1} P_2 \vfi_2 \big) = 0.
\]
By the proof of Lemma \ref{lbl4001}, $P_1 \ot \vfi_2^{-1} [P_2, \vfi_2]$ belongs to the definition ideal of the trace $\tau_1 \ot \tr_2$ and
\[
(\tau_1 \ot \tr_2) \big( P_1 \ot \vfi_2^{-1} [P_2, \vfi_2] \big) = \tau_1(P_1) \tr_2 \big( \vfi_2^{-1} [P_2, \vfi_2] \big) = \om_2(\vfi_2),
\]
since $\tau_1(P_1) = 1$. Similarly, we can show that
\[
(\tr_1 \ot \tau_2) \big( \vfi_1^{-1}[P_1, \vfi_1] \ot \vfi_2^{-1} P_2 \vfi_2 \big) = \om_1(\vfi_1).
\]

Hence, $\om(\vfi_1 \ot \vfi_2) = \theta_1 \om_1(\vfi_1) + \theta_2 \om_2(\vfi_2)$.
\end{prf}

By induction, we can generalize this result to the tensor product of any finite number of $\cstr$-algebras under similar conditions.

\begin{thm} \label{lbl4002a}
Let $n$ be a positive integer and let $\OM_i = (\CL_i, F_i, \tr_i)$ be an indicial triple for the $\cstr$-algebra $C_i$, for every $i = 1, \dots, n$. Suppose that there exists a character $\tau_i$ on $\CL_i$ such that $\tau_i(P_i) \neq 0$, where $P_i = (F_i + 1)/2$, for every $i = 1, \dots, n$. Then the triple $\OM = (\CL, F, \tr)$ is an indicial triple for $C = C_1 \ot \cdots \ot C_n$, where $\CL = \CL_1 \ot \cdots \ot \CL_n$, $F = 2(P_1 \ot \cdots \ot P_n) - 1$ and the trace $\tr: \CL^+ \to [0, +\infty]$ is defined by the relation
\[
\tr(\vfi) = \sum_{i = 1}^n \theta_i (\tau_1 \ot \cdots \ot \tau_{i - 1} \ot \tr_i \ot \tau_{i + 1} \ot \cdots \ot \tau_n) (\vfi) \quad (\vfi \in \CL^+),
\]%\tr(\vfi) = \theta_1 (\tr_1 \ot \tau_2 \ot \cdots \ot \tau_n) (\vfi) + \cdots + \theta_n (\tau_1 \ot \cdots \ot \tau_{n - 1} \ot \tr_n) (\vfi) \ \ (\vfi \in \CL^+),
where $\theta_1, \dots, \theta_n$ are fixed positive real numbers. Moreover,
\[
 \om(\vfi_1 \ot \cdots \ot \vfi_n) = \theta_1 \om_1(\vfi_1) + \cdots + \theta_n \om_n(\vfi_n),
\]
for all $\vfi_1 \in \Inv(C_1), \dots, \vfi_n \in \Inv(C_n)$, where $\om_i, \om$ are the topological indices related to the indicial triples $\OM_i, \OM$, respectively.
\end{thm}

\begin{rem} \label{lbl4002b}
In Theorem \ref{lbl4002} and Theorem \ref{lbl4002a}, it is obvious that the value of the index on tensor products at elementary invertible tensors is independent of the choice of the characters.
\end{rem}

The construction of the indicial triple for the $\cstr$-algebra $C = C_1 \ot C_2$ in Theorem \ref{lbl4002} (and similarly in Theorem \ref{lbl4002a}) depends on the existence of the characters $\tau_1, \tau_2$. We will use the following lemma to show that these characters do exist under fairly general conditions.

\begin{lem} \label{lbl4005a}
Let $\OM = (\CL, F, \tr)$ be an indicial triple for the commutative $\cstr$-algebra $C$ and let $P = {(F + 1)/2}$. If $\CL = \cstr(C, P)$ (the $\cstr$-algebra generated by $C$ and $P$), then $\Com(\CL) \subseteq \CK_\tr$.
\end{lem}

\begin{prf}
Let $C_\OM = \set{\vfi \in C \mid \df \vfi = 2 [P, \vfi] \in \CM_\tr}$. If $\vfi \in C$, then, by the density of $C_\OM$, $[P, \vfi] \in \CK_\tr$. Denote by $\CS$ the set of all elements of $\CL$ of the form
\[
a = \vfi_0 P \vfi_1 P \vfi_2 \cdots \vfi_{n-1} P \vfi_n,
\]
where $n$ is a non-negative integer and $\vfi_i \in C$, for all $0 \leq i \leq n$. Hence, for such an element $a \in \CS$, we have%
{\setlength\arraycolsep{2pt}
\begin{eqnarray*}
[a, P] &=& [\vfi_0, P] P \vfi_1 P \vfi_2 \cdots \vfi_{n-1} P
\vfi_n + \vfi_0 P [\vfi_1, P] P \vfi_2 \cdots \vfi_{n-1} P \vfi_n +\\
        && \cdots + \vfi_0 P \vfi_1 P \vfi_2 \cdots \vfi_{n-1} P [\vfi_n, P].
\end{eqnarray*}}%
Moreover, if $\vfi \in C$, then%
{\setlength\arraycolsep{2pt}
\begin{eqnarray*}
[a, \vfi] &=& \vfi_0 [P, \vfi] \vfi_1 P \vfi_2 \cdots \vfi_{n-1} P
\vfi_n + \vfi_0 P \vfi_1 [P, \vfi] \vfi_2 \cdots \vfi_{n-1} P \vfi_n + \\
          && \cdots + \vfi_0 P \vfi_1 P \vfi_2 \cdots \vfi_{n-1} [P, \vfi] \vfi_n.
\end{eqnarray*}}%
Since $\CK_\tr$ is an ideal, we have $[a, P], [a, \vfi] \in
\CK_\tr$ for any $a \in \CS, \vfi \in C$. Note that we need the commutativity of $C$ to show that the latter equality holds. Now, if $b \in \CS$, then%
{\setlength\arraycolsep{2pt}
\begin{eqnarray*}
[a, b] &=& [\vfi_0, b] P \vfi_1 P \vfi_2 \cdots \vfi_{n-1}P
       \vfi_n + \vfi_0 [P, b] \vfi_1 P \vfi_2 \cdots \vfi_{n-1} P \vfi_n +\\
       && \vfi_0 P [\vfi_1, b] P \psi_2 \cdots \vfi_{n-1} P \vfi_n  + \cdots
       + \vfi_0 P \vfi_1 P \vfi_2 \cdots \vfi_{n-1} P [\vfi_n, b].
\end{eqnarray*}}%
Thus, $[a, b] \in \CK_\tr$, for all $a, b \in \CS$.
Denote by $\CD$ the set of all the sums $\sum_{i = 1}^m a_i$, where $a_i \in \CS, \ 1 \leq i \leq m$, for some positive integer $m$. Note that $\CD$ is a  unital $*$-subalgebra of $\CL$ and $C \cup \set P \subseteq \CD \subseteq \CL$. Hence, $\CD$ is dense in $\CL$. Let $a, b \in \CD$. Then $a = \sum_{i = 1}^m a_i$ and $b = \sum_{j = 1}^n b_j$, where $a_i, b_j \in \CS$ for $1 \leq i \leq m$ and $1 \leq j \leq n$, for some positive integers $m, n$. Thus, $[a, b] = \sum_{i, j = 1}^{m, n} [a_i, b_j] \in \CK_\tr$. Hence, by the density of $\CD$, we can easily show that $[a, b] \in \CK_\tr$, for any pair of elements $a, b$ in $\CL$. This implies that the closed commutator ideal of $\CL$ is a subset of $\CK_\tr$.
\end{prf}

Theorem \ref{lbl4005} below shows that the characters needed for the construction of the trace in Theorem \ref{lbl4002} and Theorem \ref{lbl4002a} do exist under the conditions of Lemma \ref{lbl4005a}.

\begin{thm} \label{lbl4005}
Let $\OM = (\CL, F, \tr)$ be an indicial triple for the commutative $\cstr$-algebra $C$ such that $\CL = \cstr(C, P)$, where $P = {(F + 1)/2}$. Then there exists a character $\tau$ on $\CL$ such that $\tau (P) \neq 0$.
\end{thm}

\begin{prf}
By Lemma \ref{lbl4005a}, $\Com(\CL) \subseteq \CK_\tr$. Let $\pi : \CL \to \CL / \Com(\CL)$ be the quotient $*$-homomorphism from $\CL$ to the commutative unital $\cstr$-algebra $\CL / \Com(\CL)$. In the previous section we assumed that $P \notin \CM_\tr$ and this implies that $P \notin \CK_\tr$. Therefore, ${\pi(P) \neq 0}$. Hence, there exists a character $\tau_0$ on $\CL / \Com(\CL)$ such that $\tau_0(\pi(P)) \neq 0$. Thus, the mapping $\tau = \tau_0 \pi$ is a character on $\CL$ such that $\tau(P) \neq 0$.
\end{prf}

Note that the indicial triples in Examples 2.7, 2.8, 2.10 and 2.l1 in \cite{MUR1} satisfy the conditions of Theorem \ref{lbl4005}.
%Thus, we can define an indicial triple on the tensor product of any finite collection of their symbol $\cstr$-algebras.
The following examples show that the conditions in Theorem \ref{lbl4005} are essential for its validity.

\begin{exmp} \label{lbl4006}
Let $\OM_\IT = (\CL_\IT , F_\IT, \tr_\IT)$ be the indicial triple for the $\cstr$-algebra $C_\IT = C(\IT)$ discussed in Example 2.7 in \cite{MUR1}. By \cite[Theorem 2.9]{MUR1}, $\OM_2 = (\CL_2 , F_2, \tr_2)$ is an indicial triple for the $\cstr$-algebra $C_2 = \MTRX_2(C_\IT) = C_\IT \ot \MTRX_2(\IC)$, where $\CL_2 = \MTRX_2(\CL_\IT) = \CL_\IT \ot \MTRX_2(\IC)$, $F_2 = F_\IT \ot 1_2$, ($1_2$ is the identity $2 \by 2$ matrix) and the trace $\tr_2 : \CL_2 \to [0, +\infty]$ is defined as $\tr_2(\vfi) = \tr_\IT(\vfi_{11}) + \tr_\IT(\vfi_{22})$, for all $\vfi =
[\vfi_{ij}] \in \CL_2^+$. Note that $\cstr(C_2, P_2) = \CL_2$, where $P_2 = (F_2 + 1)/2 = P_\IT \ot 1_2$. There is no characters on the $\cstr$-algebra $\CL_2$. To see this, suppose that there exists a character $\tau$ on $\CL_2$. Then the mapping $\MTRX_2(\IC) \to \IC, \ x \mapsto \tau(1 \ot x)$ is a character on $\MTRX_2(\IC)$ which is impossible. Hence, there are no characters on the $\cstr$-algebra $\CL_2$. Therefore, we can not drop the commutativity condition in Theorem~\ref{lbl4005}.
\end{exmp}

\begin{exmp} \label{lbl4007}
Let $\OM_2 = (\CL_2 , F_2, \tr_2)$ be as in Example \ref{lbl4006}. Define the unital $\cstr$-subalgebra $C' = \set{\vfi \ot 1_2 \mid \vfi \in C_\IT}$ of $\CL$. It is easy to show that $\OM_2$ is an indicial triple for $C'$. Note that the $\cstr$-algebra $C'$ is commutative but there are no characters on $\CL_2$ as we have seen in Example \ref{lbl4006}. Hence, the condition $\CL = \cstr(C, P)$ in Theorem \ref{lbl4005} is essential.
\end{exmp}

\begin{rem} \label{lbl4008}
Let $\OM = (\CL, F, \tr)$ be an indicial triple for the commutative $\cstr$-algebra $C$ and let $P = (F + 1)/2$. If $\CL \neq \cstr(C, P)$, then we can define the indicial triple $\OM' = (\CL', F, \tr')$, where $\CL' = \cstr(C, P)$ and $\tr'$ is the restriction of the trace $\tr$ to $\CL'$. Both of these indicial triples has the same Toeplitz elements and the same Toeplitz algebra. Moreover, the topological index and the Fredholm index will remain unchanged too. Thus, we can apply Theorem \ref{lbl4005} by making this change in the definition of $\OM$ to obtain the character needed for the construction in Theorem \ref{lbl4002}.
\end{rem}

%___________________________________________________________________________________________________________________________
%

\section{An Index Theorem for Toeplitz Operators on the Quarter-Plane}
\label{sec3}
In this section we return to the case of Toeplitz operators on the quarter-plane discussed in the previous section to prove an index theorem for them. Let $\OM_\IT = (\CL_\IT, F_\IT, \tr_\IT)$ be the indicial triple introduced in \cite[Example 2.7]{MUR1} for the $\cstr$-algebra $C_\IT = C(\IT)$ and let $\om_\IT$ be the unique topological index related to $\OM_\IT$. By Theorem \ref{lbl4005}, there exists a character $\tau$ on $\CL_\IT$ such that $\tau(P_\IT) \neq 0$.

Let $C_{\IT^2} = C(\IT^2)$. Then $C_{\IT^2} = C_\IT \ot C_\IT$. As in the previous section, we set $\CL_{\IT^2} = \CL_\IT \ot \CL_\IT$. Since $F_{\IT^2} = 2P_{\IT^2} - 1$, where $P_{\IT^2} = P_\IT \ot P_\IT \in \CL_{\IT^2}$, then by Theorem \ref{lbl4002}, the triple $\OM_1 = (\CL_{\IT^2}, F_{\IT^2}, \tr_1)$ is an indicial triple for the $\cstr$-algebra $C_{\IT^2}$, where $\tr_1$ is the lower semicontinuous trace defined as
\[
\tr_1 : \CL_{\IT^2}^+ \To [0, \infty], \quad \tr_1(a) = (\tau \ot \tr_\IT)(a) + \theta (\tr_\IT \ot \tau)(a) \quad (a \in
\CL_{\IT^2}^+),
\]
where $\theta$ is a fixed positive irrational number (we will see later in Corollary \ref{lbl4013} the reason for making such a choice for $\theta$). To prove Theorem \ref{lbl4012} below we will redefine the trace $\tr_1$. We need the following lemma to do this.

\begin{lem} \label{lbl4009a}
In the $\cstr$-algebra $\CL_\IT$, we have $\CK_{\tr_\IT} = \CPT{L^2(\IT)} = \Com(\CL_\IT)$.
\end{lem}

\begin{prf}
By the proof of Theorem 3.3 in \cite{MUR1}, we have $\CK_{\tr_\IT} = \CPT{L^2(\IT)}$.

By Lemma \ref{lbl4005a}, $\Com(\CL_\IT) \subseteq \CK_{\tr_\IT}$. The ideal $\Com(\CL_\IT)$ contains the rank one operator $[P_\IT, z]$, therefore it contains all rank one operators. Hence, by \cite[Theorem 2.4.6]{MUR3}, $\Com(\CL_\IT)$ contains all finite rank operators and this implies that $\CPT{L^2(\IT)} \subseteq \Com(\CL_\IT)$. Hence, $ \Com(\CL_\IT) = \CK_{\tr_\IT}$.
\end{prf}

It is easily shown that the $\cstr$-algebras $\CL_\IT \ot \CK_{\tr_\IT}, \CK_{\tr_\IT} \ot \CL_\IT$ are ideals in
$\CL_{\IT^2}$. Hence, the $*$-subalgebra
\[
I = \CL_\IT \ot \CK_{\tr_\IT} + \CK_{\tr_\IT} \ot \CL_\IT = \set{a + b \mid a \in \CL_\IT \ot \CK_{\tr_\IT}, b \in \CK_{\tr_\IT} \ot \CL_\IT}
\]
is an ideal in $\CL_{\IT^2}$. Moreover, by \cite[Theorem 3.1.7]{MUR3}, $I$ is a closed ideal in $\CL_{\IT^2}$. Let $a \in \CL_\IT, b \in \CM_{\tr_\IT}$. Then $|b| \in \CM_{\tr_\IT}^+$, by \cite[Proposition A1(a)]{PHRA1}. By Lemma \ref{lbl4009a}, $|b| \in \CK_{\tr_\IT} = \Com(\CL_\IT)$. Thus, $\tau(|b|) = 0$. Hence, it follows that $|a \ot b| = |a| \ot |b| \in \CM_{\tr_1}^+$. Again by \cite[Proposition A1(a)]{PHRA1}, $a \ot b \in \CM_{\tr_1}$. Thus, $\CL_\IT \ot \CK_{\tr_\IT} \subseteq \CK_{\tr_1}$. Similarly, we can show that $\CK_{\tr_\IT} \ot \CL_\IT \subseteq \CK_{\tr_1}$. Hence, $I \subseteq \CK_{\tr_1}$. The mapping $\tr_{\IT^2} : \CL_{\IT^2}^+ \to [0, +\infty]$ is a lower semicontinuous trace on $\CL_{\IT^2}$, where
\[
\tr_{\IT^2}(a) = \left\{ \begin{array}{ll}
\tr_1(a), & \textrm{if $a \in I^+$,}\\
+\infty, & \textrm{if $a \in \CL_{\IT^2}^+ \backslash I^+$.}\\
\end{array} \right.
\]
Moreover, $\CK_{\tr_{\IT^2}} = I$. Let $\OM_{\IT^2}$ be the triple $(\CL_{\IT^2}, F_{\IT^2}, \tr_{\IT^2})$. If $\vfi_1, \vfi_2 \in C_{\OM_\IT}$, then
\[
\df (\vfi_1 \ot \vfi_2) = 2 ( [P_\IT, \vfi_1] \ot P_\IT \vfi_2 + \vfi_1 P_\IT \ot [P_\IT, \vfi_2]) \in \CM_{\tr_1} \cap I =
\CM_{\tr_{\IT^2}}.
\]
Hence, the triple $\OM_{\IT^2} = (\CL_{\IT^2}, F_{\IT^2}, \tr_{\IT^2})$ is an indicial triple for $C_{\IT^2} = C(\IT^2)$. Denote by $\om_{\IT^2}$ the unique topological index related to $\OM_{\IT^2}$.

Let $\vfi \in \Inv(C_{\IT^2})$. By a theorem of van Kampen \cite{VKM1}, there exists $\psi \in C_{\IT^2}$ and a unique
character $\chi$ on the group $\IT^2$ such that $\vfi = \chi e^\psi$. By \cite[Theorem 2.2.2]{RUD1}, $\chi = z^m \ot z^n$, for a unique pair of integers $m, n \in \IZ$. Thus, we have $\vfi = (z^m \ot z^n) e^\psi$. To calculate the topological index for $\vfi$, first note that $\om_{\IT^2} ( \vfi ) = \om_{\IT^2} ( ( z^m \ot z^n ) e^\psi ) = \om_{\IT^2} ( z^m \ot z^n ) + \om_{\IT^2} ( e^\psi ) = \om_{\IT^2} ( z^m \ot z^n )$, since $\om_{\IT^2} ( e^\psi ) = 0$. Let $\om_1$ be the unique topological index on $C_{\IT^2}$ related to the indicial triple $\OM_1$ defined above. Since $\df (z^m \ot z^n) \in \CM_{\tr_{\IT^2}}$, we have {\setlength\arraycolsep{2pt}
\begin{eqnarray*}
   \om_{\IT^2}(z^m \ot z^n)&=& \frac12 \tr_{\IT^2}((z^m \ot z^n)^{-1} \df (z^m \ot z^n)) \\
                           &=& \frac12 \tr_1((z^m \ot z^n)^{-1} \df (z^m \ot z^n)) \\
                           &=& \om_1(z^m \ot z^n) \\
                           &=& \theta \om_\IT(z^m) + \om_\IT(z^n) \quad \textrm{(by Theorem \ref{lbl4002})} \\
                           &=& \theta m \om_\IT(z) + n \om_\IT(z) = \theta m + n \quad \textrm{(by \cite[Example 2.7]{MUR1})}.
\end{eqnarray*}}%
Hence, we have proved the following theorem.

\begin{thm} \label{lbl4010}
If $\vfi \in \Inv(C_{\IT^2})$, then the topological index of $\vfi$ is $\om_{\IT^2} (\vfi) = \theta m + n$, where $m, n$ are the unique pair of integers such that $\vfi = (z^m \ot z^n) e^\psi$, for some $\psi \in C_{\IT^2}$.
\end{thm}

\begin{rem} \label{lbl4011a}
It follows immediately from Theorem \ref{lbl4010} that the value of the topological index $\om_{\IT^2}$ at any $\vfi \in \Inv(C_{\IT^2})$ is independent of the choice of the character $\tau$. %Note that the topological index $\om_{\IT^2}$ maps the set $\Inv(C_{\IT^2})$ onto a subgroup of $\IR$ isomorphic to the group $\IZ^2$.
\end{rem}

Denote by $\CA_{\IT^2}$ the Toeplitz algebra associated with the indicial triple $\OM_{\IT^2}$ and denote by $\Tr_{\IT^2}$ the restriction of the trace $\tr_{\IT^2}$ to $\CA_{\IT^2}$. The $\cstr$-algebra $\CA_{\IT^2}$ is $*$-isomorphic to the Toeplitz algebra $\CT_{\IT^2}$ generated by Toeplitz operators on the quarter-plane defined in the introduction. Thus, we can identify these two $\cstr$-algebras. By \cite[Theorem 3.1]{MUR1}, the following theorem is an immediate corollary of Theorem \ref{lbl4010}.

\begin{thm} \label{lbl4011}
If $\vfi \in \Inv(C_{\IT^2})$, then the Toeplitz element $T_\vfi$ associated to the indicial triple $\OM_{\IT^2}$ for $C_{\IT^2}$ is an $\OM_{\IT^2}$-Fredholm element and the Fredholm index of $T_\vfi$ relative to $\Tr_{\IT^2}$ is
\[
\ind_{\OM_{\IT^2}} (T_\vfi) = -\theta m - n,
\]
where $m, n$ are the unique integers such that $\vfi = (z^m \ot z^n) e^\psi$, for some $\psi \in C_{\IT^2}$.
\end{thm}

The following theorem gives a characterization of $\OM_{\IT^2}$-Fredholm Toeplitz elements. Before proving it, we give some properties of the Toeplitz algebra $\CA_{\IT^2}$. By \cite[Proposition 1]{DOHO1}, $\CA_{\IT^2} = \CA_\IT \ot \CA_\IT$, where $\CA_\IT$ is the Toeplitz algebra associated to the indicial triple $\OM_\IT$. By \cite[p. 205]{DOHO1}, there is a unital $*$-homomorphism $\pi_{\IT^2} : \CA_{\IT^2} \to C_{\IT^2}$ such that $\pi_{\IT^2}(T_\vfi) = \vfi$, for all $\vfi \in C_{\IT^2}$, and $\ker(\pi_{\IT^2}) = \Com(\CA_{\IT^2})$, where $\Com(\CA_{\IT^2})$ is the closed commutator ideal of $\CA_{\IT^2}$. Moreover, by \cite[p. 207]{DOHO1},
\[
\Com(\CA_{\IT^2}) = \CA_\IT \ot \Com(\CA_\IT) + \Com(\CA_\IT) \ot \CA_\IT.
\]
From the theory of classical Toeplitz operators, $\Com(\CA_\IT) = P_\IT \CPT{L^2(\IT)} P_\IT = \CK_{\Tr_\IT}$, where $\Tr_\IT$ is the restriction of the trace $\tr_\IT$ to $\CA_\IT$.

\begin{thm} \label{lbl4012}
Let $\vfi \in C_{\IT^2}$. Then the Toeplitz element $T_\vfi$ is an $\OM_{\IT^2}$-Fredholm element if, and only if, $\vfi \in \Inv(C_{\IT^2})$.
\end{thm}

\begin{prf}
First we show that $\CK_{\Tr_{\IT^2}} \subseteq  \Com(\CA_{\IT^2})$. If $S \in \CK_{\Tr_{\IT^2}}$, then $S \in \CK_{\tr_{\IT^2}}$. We can write $S = a + b$, where $a \in \CL_\IT \ot \CK_{\tr_\IT}, b \in \CK_{\tr_\IT} \ot \CL_\IT$, since $\CK_{\tr_{\IT^2}} = \CL_\IT \ot \CK_{\tr_\IT} + \CK_{\tr_\IT} \ot \CL_\IT$.

By density, we can choose two sequences $(a_n)_{n \geq 1}$ in $\CL_\IT \odot \CK_{\tr_\IT}$ and $(b_n)_{n \geq 1}$ in
$\CK_{\tr_\IT} \odot \CL_\IT$, such that the sequence $a_n \to a$ and $b_n \to b$. It can be shown that $P_\IT u P_\IT \in \CA_\IT$, whenever $u \in \CL_\IT$. Similarly, $P_\IT v P_\IT \in \CK_{\Tr_\IT}$, whenever $v \in \CK_{\tr_\IT}$. Thus, $P_{\IT^2} a_n P_{\IT^2} \in \CA_\IT \odot \CK_{\Tr_\IT} \subseteq \CA_\IT \ot \CK_{\Tr_\IT}$, for all $n \geq 1$. Hence,
\[
P_{\IT^2} a P_{\IT^2} \in \CA_\IT \ot \Com(\CA_\IT) \subseteq
\Com(\CA_{\IT^2}).
\]
Similarly, we can show that $P_{\IT^2} b P_{\IT^2} \in \Com(\CA_{\IT^2})$. Therefore, $S = P_{\IT^2} S P_{\IT^2} \in
\Com(\CA_{\IT^2})$. Thus, $\CK_{\Tr_{\IT^2}} \subseteq \Com(\CA_{\IT^2}) = \ker(\pi_{\IT^2})$.

Now if $T_\vfi$ is an $\OM_{\IT^2}$-Fredholm element, then there is an element $S \in \CA_{\IT^2}$ such that both of the elements $P_{\IT^2} - T_\vfi S, P_{\IT^2} - S T_\vfi$ belong to the ideal $\CK_{\Tr_{\IT^2}}$. It follows immediately that $\vfi \pi_{\IT^2}(S) = \pi_{\IT^2}(S) \vfi = 1$. Hence, $\vfi$ is an invertible function. The opposite direction follows from Theorem \ref{lbl4011}.
\end{prf}

\begin{cor} \label{lbl4013}
Let $\vfi \in C_{\IT^2}$. Then the Toeplitz element $T_\vfi$ is an $\OM_{\IT^2}$-Fredholm element with $\ind_{\OM_{\IT^2}} (T_\vfi) = 0$ if, and only if, $\vfi = e^\psi$, for some $\psi \in C_{\IT^2}$.
\end{cor}

\begin{prf}
Let $\vfi \in C_{\IT^2}$ and  $T_\vfi$ is an $\OM_{\IT^2}$-Fredholm element with zero $\OM_{\IT^2}$-Fredholm index. By theorem \ref{lbl4012}, $\vfi \in \Inv(C_{\IT^2})$. Thus, $\vfi = (z^m \ot z^n) e^\psi$, for a unique pair of integers $m, n$ and for some $\psi \in C_{\IT^2}$. By Theorem \ref{lbl4011}, $-\theta m - n = \ind_\OM (T_\vfi) = 0$. Hence, $m = n = 0$, since $\theta$ is irrational. Thus, $\vfi = e^\psi$. The reverse implication is clear.
\end{prf}

\begin{rem} \label{lbl4014}
From this corollary we can see the reason for the choice of $\theta$ as an irrational positive number.
\end{rem}

In the theory of Toeplitz operators on the unit circle, a Toeplitz operator is invertible if, and only if, it is a Fredholm operator and its Fredholm index is zero \cite[Theorem 3.5.15]{MUR3}. The situation is different in the case of Toeplitz operators on the quarter-plane. Every invertible Toeplitz operator on the quarter-plane is an $\OM_{\IT^2}$-Fredholm element and its $\OM_{\IT^2}$-index is zero. However, in \cite[p. 210]{DOHO1} there is an example of a function $\vfi \in \Inv(C_{\IT^2})$ such that $T_\vfi$ is a Fredholm operator in the classical sense and its Fredholm index is zero but $T_\vfi$ is not invertible. Therefore, by \cite[Corollary of Theorem 1]{DOHO1}, $\vfi = e^\psi$, for some $\psi \in C_{\IT^2}$. Thus, by Theorem \ref{lbl4011}, $T_\vfi$ is an $\OM_{\IT^2}$-Fredholm element and $\ind_{\OM_{\IT^2}} (T_\vfi) = 0$. Hence, there is an $\OM_{\IT^2}$-Fredholm Toeplitz element such that its $\Tr_{\IT^2}$-index is zero but it is not invertible. However, we can show that an $\OM_{\IT^2}$-Fredholm Toeplitz element with zero $\Tr_{\IT^2}$-index has a property weaker than invertibility. We give this property in the following proposition.

\begin{prop} \label{lbl4016}
Let $\vfi \in C_{\IT^2}$ and let $T_\vfi$ be an $\OM_{\IT^2}$-Fredholm element such that  $\ind_{\OM_{\IT^2}}(T_\vfi) = 0$. Then there exists an invertible element $S \in \CA_{\IT^2}$ such that $T_\vfi - S \in \CK_{\Tr_{\IT^2}}$.
\end{prop}

\begin{prf}
By Corollary \ref{lbl4013}, $\vfi = e^\psi$, for some $\psi \in C_{\IT^2}$. Let $\pi_{\IT^2}: \CA_{\IT^2} \to C_{\IT^2}$ be the unique $*$-homomorphism such that $\pi_{\IT^2}(T_\vfi) = \vfi$, for all $\vfi \in C_{\IT^2}$. Then $\pi_{\IT^2}(T_\vfi) = e^\psi$. If we choose $S_0 \in \CA_{\IT^2}$ such that $\pi_{\IT^2}(S_0) = \psi$, then $\pi_{\IT^2}(e^{S_0}) = e^{\pi(S_0)} = e^\psi = \pi_{\IT^2}(T_\vfi)$. Thus, $T_\vfi - e^{S_0} \in \CK_{\Tr_{\IT^2}}$. Thus, we have proved our result since $e^{S_0}$ is an invertible element in $\CA_{\IT^2}$.
\end{prf}

\begin{rem}
We have defined Toeplitz operators on $H^2(\IT^N)$ in the introduction. Using Theorem \ref{lbl4002a}, we can establish an index theorem for these operators similar to the one we have proved for Toeplitz operators on the quarter-plane. The construction and proofs are similar.
\end{rem}

Now we come to the extension of the previous results to Toeplitz operators with matrix symbols acting on the Hardy space $H_N^2(\IT^2) = \bigoplus_{k = 1}^N H^2(\IT^2)$ defined in the introduction. Let $\OM_{\IT^2}$ be the indicial triple for $C(\IT^2)$ defined above. Define $C_N = \MTRX_N(C(\IT^2)) = C(\IT^2) \ot \MTRX_N(\IC)$, $\CL_N = \MTRX_N(\CL_{\IT^2}) = \CL_{\IT^2} \ot \MTRX_N(\IC)$ and define the trace $\tr_N: \CL_N^+ \to [0, +\infty]$ as $\tr_N([a_{ij}]) = \sum_{i = 1}^N \tr_{\IT^2}(a_{ii})$, where $[a_{ij}] \in \CL_N^+$. Let $F_N = F_{\IT^2} \ot 1_N$, where $1_N$ is the unit element in the algebra $\MTRX_N(\IC)$. Note that $P_N = (F_N + 1)/2$. Then, by \cite[Theorem 2.9]{MUR1}, $\OM_N = (\CL_N, F_N, \tr_N)$ is an indicial triple for $C_N$. Moreover, if $\Phi \in \Inv(C_N)$, then $\det \Phi$ is a non-vanishing function on $\IT^2$ and $\om_N(\Phi) = \om_{\IT^2} (\det \Phi) = \theta m + n$, where $m, n$ are the unique integers such that $\det \Phi = (z^m \ot z^n) e^\psi$, for some $\psi \in C_{\IT^2}$ and $\om_N$ is the topological index related to $\OM_N$. Let $\CA_N$ be the Toeplitz algebra related to $\OM_N$ and let $\Tr_N$ be the restriction of the trace $\tr_N$ to $\CA_N$. Now we prove the index theorem for Toeplitz operators with symbols belonging to $C_N$.

\begin{thm}
Let $\Phi \in C_N$. Then $T_\Phi$ is an $\OM_N$-Fredholm operator if and only if $\det \Phi$ is non-vanishing on $\IT^2$ and in this case we have
\[
\ind_N(T_\Phi) = -\theta m - n,
\]
where $m, n$ are the unique integers such that $\det \Phi = (z^m \ot z^n) e^\psi$, for some $\psi \in C_{\IT^2}$ and $\ind_N$ is the $\Tr_N$-index mapping.
\end{thm}

\begin{prf}
We define the unital $*$-homomorphism $\pi_N: \CA_N \to C_N$, such that $\pi_N(S) = [\pi_{\IT^2}(s_{ij})]$, where $S = [s_{ij}] \in \CA_N$. That is, $\pi_N$ is the inflation homomorphism of $\pi_{\IT^2}$. It is easy to show that $\pi_N(T_\Phi) = \Phi$, for all $\Phi \in C_N$ and $\ker(\pi_N) = \MTRX_N(\ker(\pi_{\IT^2}))$. By some matrix manipulations, we can show that $\CM_{\Tr_N} \subseteq \MTRX_N(\CM_{\Tr_{\IT^2}})$. Hence, by the proof of Theorem \ref{lbl4012}, $\CK_{\Tr_N} \subseteq \MTRX_N(\CK_{\Tr_{\IT^2}}) \subseteq \MTRX_N(\ker(\pi_{\IT^2}))$.

Let $\Phi \in C_N$ such that $T_\Phi$ is an $\OM_N$-Fredholm operator. Then there exists $S \in \CA_N$ such that $P_N - S T_\Phi, P_N - T_\Phi S \in \CK_{\Tr_N}$. Hence, $1 - \pi_N(S) \Phi = 1 - \Phi \pi_N(S) = 0$. Thus, it follows immediately that $\det \Phi$ is a non-vanishing function on $\IT^2$. The rest of the assertions follow immediately from \cite[Theorem 3.1]{MUR1} and the paragraph preceding this theorem.
\end{prf}

%___________________________________________________________________________________________________________________________
%

\bibliographystyle{amsplain}

\end{document}